\documentclass[a4paper,english,12pt]{article}
\usepackage[T1]{fontenc}
\usepackage[latin1]{inputenc}
\usepackage{babel}
\usepackage{amsmath}
\usepackage{amssymb}
\usepackage{makeidx}
\usepackage{fancybox}
\usepackage{amsfonts}
\usepackage{latexsym}
\usepackage{cite}
\usepackage[numbers]{natbib}

\usepackage{anysize}
\usepackage{amsthm}
\usepackage[numbers]{natbib}
\usepackage[pdftex]{graphicx}
\usepackage[pdftex]{color}

\usepackage[all]{xy}

\newtheorem{thm}{Theorem}[section]
 
 \newtheorem{lem}[thm]{Lemma}
 \newtheorem{prop}[thm]{Proposition}
 \newtheorem{defn}[thm]{Definition}

\newcommand*{\m}[1]{\underline{#1}}

\begin{document}
\title{Hypermonogenic solutions and plane waves of the Dirac operator in $\mathbb{R}^p\times\mathbb{R}^q$}
\author{Al\'{i} Guzm\'{a}n Ad\'{a}n, Heikki Orelma, Franciscus Sommen}
\maketitle

\begin{abstract}
In this paper we first define hypermonogenic solutions of the Dirac operator in $\mathbb{R}^p\times\mathbb{R}^q$ and study some basic properties, e.g., obtaining a Cauchy integral formula in the unit hemisphere. Hypermonogenic solutions form a natural function class in classical Clifford analysis. After that, we define the corresponding hypermonogenic plane wave solutions and deduce explicit methods to compute these functions.
\end{abstract}


%\tableofcontents

\section{Introduction}
Clifford analysis is nowadays a well established generalization of the classical complex analysis in higher dimensions. It is a function theory for solutions of the Dirac operator $\partial_{\m{x}}=\sum_{j=1}^me_j\partial_{x_j}$, where $\{ e_1,...,e_m\}$ produces the Clifford algebra $\mathbb{R}_m$ by the defining relations $e_ie_j+e_je_i=-2\delta_{ij}$. One of the key features of the Dirac operator is its $SO(m)$ rotation invariance. The effect of this feature may be seen in every part of the theory. This inspired also look for more general operators, acting on Clifford algebra valued functions, but with a ''weaker''  symmetry (for example a subgroup of $SO(m)$). On example of these is the so called modified Dirac operator defined by Heinz Leutwiler and Sirkka-Liisa Eriksson. Their idea is to add an extra generator $e$ and then generate the Clifford algebra $\mathbb{R}_{m+1}$ by the defining relations $e_ie_j+e_je_i=-2\delta_{ij}$, $e^2=-1$ and $ee_j=-e_je$. They obtained the operator, called the modified Dirac operator
\[
M=e\partial_r+\partial_{\m{y}}+\frac{p-1}{r}e\cdot,
\]
where $e\cdot$ denoted the interior multiplication and $p$ is a natural number. This operator is $SO(m)$ invariant only with respect to $e$-axis, not anymore to the whole space. In this paper we deduce a correspondence between the set of the null solutions of the modified Dirac operator $M$ and the solutions of the biaxial Dirac operator $\partial_{\m{x}}+\partial_{\m{y}}$. In particular, for $f=A+eB$ we obtain that
\[
\partial_{\m{x}}(A+eB)+\partial_{\m{y}}(A+eB)=e\partial_r(A+eB)+\partial_{\m{y}}(A+eB)+\frac{p-1}{r}e\cdot (A+eB),
\]
where $A$ and $B$ are radial functions with respect to $\m{x}$ taking values in $\mathbb{R}_m$ and $e=\frac{\m{x}}{r}$, where $r=|\m{x}|$.  We will call these functions hypermonogenic solutions of the Dirac operator; they are special biaxial monogenic functions. \\
\\
In the first Section 3 we study hypermonogenic functions as special biaxial solutions of the Dirac equation and we present the corresponding Cauchy-Kovalevskaya extension formula. In the next section we obtain a Cauchy formula for hypermonogenic functions in the upper half-ball, resulting from Cauchy's formula for biaxial monogenics in the unit ball. In the next section we construct a family of plane wave hypermonogenic solutions by applying Funk-Hecke's formula to an integral of plane wave monogenic function. Examples are given in the last section, in particular using the monogenic Fourier kernel $\exp(\langle \m{x},\m{t}\rangle +i\langle \m{y},\m{s}\rangle)(\m{t}+i\m{s})$. This example is expressed in terms of Bessel functions and it also can be obtained by applying the Cauchy-Kovalevskaya extension formula. This is the starting point for Fourier analysis and Radon transforms for hypermonogenic solutions.

\section{Preliminaries}
Let us now recall some standard facts from theory of Clifford algebra and analysis. We denote by $\mathbb{R}_m$ the real Clifford algebra generated by the symbols $\{e_1,...,e_m\}$ satisfying the defining relations
$
e_ie_j+e_je_i=-2\delta_{ij},
$
for $i,j=1,...,m$, where $\delta_{ij}$ is the well known Kronecker symbol. We may consider the Euclidean vector space $\mathbb{R}^m$ as a natural embedded subspace of $\mathbb{R}_m$ with the embedding $(x_1,...,x_m)\mapsto \m{x}:=\sum_{j=1}^m x_je_j$. Similarly, we define the complex Clifford algebra by $\mathbb{C}_m=\mathbb{R}_m\otimes \mathbb{C}$  and  consider $\mathbb{C}^m$ as its subspace. If we define the increasing $k$-list by $A=\{ (i_1,...,i_k): 0\le i_l\le m,\ i_1+...+i_m=k,\ i_1<...<i_k\}$ and put $|A|=k$ we may write basis elements of a Clifford algebra by $e_A:=e_{i_1}...e_{i_k}$. Then an element $a\in \mathbb{R}_m$ (or $a\in \mathbb{C}_m$) will be expressed as a sum $a=\sum_{k=0}^n\sum_{|A|=k} a_Ae_A$ where $a_A\in\mathbb{R}$ (or $a_A\in\mathbb{C}$). In $\mathbb{R}^m$ we define the euclidean inner product by $\langle \m{x},\m{y}\rangle=\sum_{j=1}^mx_jy_j$ and the hermitean inner product in $\mathbb{C}^m$  by $\langle \m{x},\m{y}\rangle=\sum_{j=1}^mx_jy^c_j$, where $y^c_j$ is the complex conjugation of ${y}_j$. An inner product induces the norm $|\m{x}|=\langle \m{x},\m{x}\rangle^{\frac{1}{2}}$.  If $\m{x},\m{y}\in \mathbb{R}_m$ are vectors, the product of these vectors may be represented as
\[
\m{x}\m{y}=\frac{1}{2}(\m{x}\m{y}+\m{y}\m{x})+\frac{1}{2}(\m{x}\m{y}-\m{y}\m{x}),
\]
it is easy to see, that the first term is a real number and $\langle \m{x},\m{y}\rangle=-\frac{1}{2}(\m{x}\m{y}+\m{y}\m{x})$. The last term is a linear combination of basis elements $e_ie_j$, where $i<j$, and denoted by $\m{x}\wedge\m{y}=\frac{1}{2}(\m{x}\m{y}-\m{y}\m{x})$. These operators may be generalized as follows. In the Clifford algebra $\mathbb{R}_m$ we may define the subspaces of $k$-multivectors as
$
\mathbb{R}_m^{(k)}=\{a_k= \sum_{|A|=k}x_Ae_A : x_A\in\mathbb{R}\}.
$
Then we may split the multiplication of a vector and of an $k$-multivector as
\[
\m{x}a_k=\m{x}\cdot a_k+\m{x}\wedge a_k,
\]
where the interior and exterior multiplication are 
\begin{align*}
\m{x}\cdot a_k=\frac{1}{2}(\m{x} a_k-(-1)^ka_k\m{x}),\\
\m{x}\wedge a_k=\frac{1}{2}(\m{x} a_k+(-1)^ka_k\m{x}).
\end{align*}
Since an arbitrary element $a\in\mathbb{R}_m$ admits the multivector decomposition $
a=\sum_{k=0}^m a_k,
$
where $a_k\in \mathbb{R}_m^{(k)}$ we may define the multiplications as
\begin{align*}
\m{x}\cdot a=\sum_{k=1}^m\m{x}\cdot a_k,\\
\m{x}\wedge a=\sum_{k=0}^{m-1}\m{x}\wedge a_k.
\end{align*}
Let us now consider functions $f:\Omega\to\mathbb{R}_m$ where $\Omega$ is an open subset of $\mathbb{R}^m$. We say that these functions, which are of the form
\[
f=\sum_{k=0}^m\sum_{|A|=k}f_Ae_A,
\]
 are differentiable, integrable etc. if their component functions $f_A$ have these properties. The fundamental differential operator in theory of Clifford analysis is the so called Dirac operator, defined by
\[
\partial_{\m{x}}=e_1\partial_{x_1}+\cdots+e_m\partial_{x_m}
\]
and  its square equals the Laplacian $-\Delta_{\m{x}}$, that is, 
\[
\Delta_{\m{x}}=-\partial_{\m{x}}^2=\partial_{x_1}^2+...+\partial_{x_m}^2.
\]
A differentiable Clifford algebra valued function $f$, defined on some open set $\Omega\subset\mathbb{R}^m$, is called (left) monogenic (resp. right) if $\partial_{\m{x}}f=0$ (resp. $f\partial_{\m{x}}=0$) on $\Omega$. Monogenic functions have many nice function theoretic properties, see \cite{BDS, DSS}. On fundamental feature is the Cauchy integral formula. In this paper we need it in the following special case. Consider monogenic functions $f:\Omega\to\mathbb{R}_m$, where $\Omega\subset\mathbb{R}^m $ contains the unit ball $B_1=\{ x\in\mathbb{R}^m : |x|<1\}$. The surface of the unit ball is the unit sphere $S^{m-1}=\{ \m{x}\in\mathbb{R}^m : |\m{x}|=1\}$. 

\begin{thm}[Cauchy formula on the unit ball, \cite{BDS}]
Let $f:\Omega\to\mathbb{R}_m$ be a monogenic function and $B_1\subset\Omega$. If $\m{x}\in B_1$, then
\[
f(\m{x})=\frac{1}{\lambda_{m-1}}\int_{S^{m-1}} \frac{\m{x}-\m{\eta}}{|\m{x}-\m{\eta}|^m}d\sigma(\m{\eta})f(\m{\eta}),
\]
where $\lambda_{m-1}$ is the measure of $S^{m-1}$ and the vectorial surface element $d\sigma(\m{\eta})=\m{\eta}dS(\m{\eta})$ is the product of the outwart pointing normal unit vector $\m{\eta}$ on the sphere and the scalar surface element $dS(\m{\eta})$.
\end{thm}
On the unit sphere we will need also the following integral formula.

\begin{thm}[Funk-Hecke, \cite{H}]
Let $\m{\xi},\m{\eta}\in S^{m-1}$ and suppose $\psi(t)$ is continuous for $-1\le t\le 1$. Then for every spherical harmonic  $H_k$  of degree $k$ 
\[
\int_{S^{m-1}} \psi(\langle\m{\xi},\m{\eta}\rangle)H_k(\m{\eta})dS(\m{\eta})= \lambda_{m-1}\frac{k!}{(m-2)_k} H_k(\m{\xi})\int_{-1}^1 \psi (t)C^{\frac{m}{2}-1}_k(t)(1-t^2)^{\frac{m-3}{2}}dt,
\]
where $C_k^\ell$ is a Gegenbauer polynomial and $(a)_k$ the Pochhammer symbol.
\end{thm}

\section{Hypermonogenic solutions  in biaxial domains}
In the beginning of 90's Heinz Leutwiler together with Sirkka-Liisa Eriksson defined the so called hypermonogenic functions (see e.g. \cite{ELI,L}).
They are functions $f:\Omega\subset \mathbb{R}^{q+1}_+\to \mathbb{R}_{q+1}$ which are defined on the open subset of the upper half-space $\mathbb{R}^{q+1}_+=\{ (r,\m{y})\in\mathbb{R}\times \mathbb{R}^{q} : r>0\}$, take values in the Clifford algebra $\mathbb{R}_{q+1}$ associated to a quadratic form $(re+\m{y})^2=-r^2-|\m{y}|^2$ and belong to the kernel of the  so called modified Dirac operator 
\[
Mf:=e\partial_rf+\partial_{\m{y}}f+\frac{p-1}{r}e\cdot f,
\]
where $e\cdot f$ denotes the interior multiplication operator and $p$ a natural number. \\
\\
The theory of hypermonogenic functions, or briefly, the hypermonogenic function theory is related to a hyperbolic upper half-space model. The theory have it's own advantages, e.g., paravector power functions may be included in the kernel of the operator. Also integral formulas are studied  and a wide class of special solutions. The reader may find more information and recent research results for example in \cite{ELI, EOS, EOV,L}. \\
\\
Let us now show, how these hypermonogenic functions may be seen as a special case of monogenic functions. Our fundamental idea is to consider $r$ as a norm of some extra $p$-vector variable $\m{x}$, that is $r=|\m{x}|$. Considering only radial functions with respect to $\m{x}$, we have $\partial_{\m{x}}B=\frac{\m{x}}{r}\partial_rB$ and we may denote $e=\frac{\m{x}}{r}$. Now $\partial_{\m{x}}e=\frac{1-p}{r}$ and we obtain for radial functions $\partial_{\m{x}}(eB)=-\partial_rB-\frac{p-1}{r}B=e\partial_r(eB)+\frac{p-1}{r}e\cdot (eB)$. If we consider now functions $A$ and $B$, which are radial with respect to $\m{x}$ and take their values in a Clifford algebra $\mathbb{R}_{q}$ associated to a quadratic form $\m{y}^2=-|\m{y}|^2$, we have that $\partial_{\m{x}}(A+eB)+\partial_{\m{y}}(A+eB)=e\partial_r(A+eB)+\partial_{\m{y}}(A+eB)+\frac{p-1}{r}e\cdot (A+eB)$. This motivates us to make the following definition.

\begin{defn}[Hypermonogenic solutions] Let $f:\Omega\subset \mathbb{R}^p\times\mathbb{R}^q\to \mathbb{C}_{p+q}$ be a differentiable function. 
Hypermonogenic solutions of the Dirac operator are functions of the form
\[
f(\m{x},\m{y})=A(|\m{x}|,\m{y})+\frac{\m{x}}{|\m{x}|}B(|\m{x}|,\m{y}),
\]
defined on a biaxial domain $\Omega$ in $\mathbb{R}^p\times\mathbb{R}^q$ satisfying $(\partial_{\m{x}}+\partial_{\m{y}})f=0$.
\end{defn}
Let us now find a generalized Cauchy-Riemann system for hypermonogenic solutions.

\begin{prop}
A function $f(\m{x},\m{y})=A(|\m{x}|,\m{y})+\frac{\m{x}}{|\m{x}|}B(|\m{x}|,\m{y})$ is a hypermonogenic solution, if and only if
\begin{align*}
\partial_{\m{y}}A-\partial_{|\m{x}|}B-\frac{p-1}{|\m{x}|}B=0,\\
\partial_{\m{y}}B-\partial_{|\m{x}|}A=0.
\end{align*}
\end{prop}
Proof. Since $\partial_{\m{x}}(\frac{\m{x}}{|\m{x}|})=\frac{1-p}{|\m{x}|}$ we obtain
\begin{align*}
\partial_{\m{x}}\Big(\frac{\m{x}}{|\m{x}|}B(|\m{x}|,\m{y})\Big)
&=\partial_{\m{x}}\big(\frac{\m{x}}{|\m{x}|})B(|\m{x}|,\m{y})+\sum_{j=1}^{p}e_j\frac{\m{x}}{|\m{x}|}\partial_{x_j}B(|\m{x}|,\m{y})\\
&=\frac{1-p}{|\m{x}|}B(|\m{x}|,\m{y})+\frac{\m{x}^2}{|\m{x}|^2}\partial_{|\m{x}|}B(|\m{x}|,\m{y}).
\end{align*}
The system is easy obtain using this formula.$\square$\\
\\
In the heart of this paper are functions with the series representation
\begin{align}\label{CHoj}
f(\m{x},\m{y})&=\sum_{j=0}^\infty \m{x}^j f_j(\m{y}).
\end{align}
All monogenic functions 
of this form\footnote{One may find, that in this case $A=\sum_{j=0}^\infty(-1)^{j}|\m{x}|^{2j}f_{2j}(\m{y})$ and  $ B=\sum_{j=0}^\infty(-1)^{j}|\m{x}|^{2j+1}f_{2j+1}(\m{y})\nonumber$.} are hypermonogenic solutions. For example, in the second part of the paper, these type of solutions are crucial when we consider the so called plane wave solutions. We can also find hypermonogenic solutions of this form as follows.\\
\\
The Vekua type system in the preceding proposition allows us to construct hypermonogenic solutions using Cauchy-Kovalevskaya extension in higher codimensions, see \cite{DSS}. Hypermonogenic solutions are determined by the restriction  $A(0,\m{y})$ to the surface $\m{x}=\m{0}$, that is, starting from a given function $A(\m{y})$ we may establish its Cauchy-Kovalevskaya extension 
\[
A(|\m{x}|,\m{y})+\frac{\m{x}}{|\m{x}|}B(|\m{x}|,\m{y}).
\]
Let us now take an open set $\widetilde{\Omega}\subset \mathbb{R}^p\times\mathbb{R}^q$, which is $SO(p)$ invariant and denote $\Omega=\widetilde{\Omega}\cap \mathbb{R}^q$ assuming this to be non empty. In a neighbourhood of this set, we may construct a hypermonogenic function, when $f_{0}(\m{y})$ is given, as follows.

\begin{thm}[Generalized Cauchy-Kovalevskaya  extension  \cite{DSS}]\label{Absolut} Let $f_{0}(\m{y})$ be an analytic function in $\Omega$. Then there exists a unique sequence $\{f_j(\m{y})\}_{j=1}^\infty$ of analytic functions such that the series
\[
f(\m{x},\m{y})=\sum_{j=0}^\infty \m{x}^j f_j(\m{y})
\]
is convergent in a neigbourhood  $U\subset \mathbb{R}^p\times \mathbb{R}^q$ of $\Omega$ and its sum $f$ is a hypermonogenic solution in $U$. The function $f_{0}(\m{y})$ is determined by the relation
\[
f_{0}(\m{y})=f(\m{0},\m{y}).
\]
Furthermore, the sum $f$ is formally given by the expression
\begin{align*}
f(\m{x},\m{y})=\Gamma(p/2)\Big(\frac{|\m{x}|}{2}\sqrt{\Delta_{\m{y}}}\Big)^{-p/2}\Big(\frac{|\m{x}|}{2}\sqrt{\Delta_{\m{y}}}J_{\frac{p}{2}-1}(|\m{x}|\sqrt{\Delta_{\m{y}}})+\frac{\m{x}\partial_{\m{y}}}{2}J_{\frac{p}{2}}(|\m{x}|\sqrt{\Delta_{\m{y}}})\Big)f_{0}(\m{y}),
\end{align*}
where $J$ is a Bessel function and $\sqrt{\Delta_{\m{y}}}$ is the the formal square root of the Laplacian.
\end{thm}
Also in the paper \cite{DS} the Cauchy-Kovalevskaya problem in the biaxial framework  was studied. The aim of that paper was to find a hypermonogenic solution of the form
\[
f(\m{x},\m{y})=\sum_{j=0}^\infty c_j\m{x}^j F_j(\m{y})W(\m{y})
\]
where $F_j(\m{y})$ are weight functions and $W(\m{y})$ is a given initial function whereby $F_0(\m{y})=1$, $c_0=1$ and $\partial_{\m{x}}(c_j\m{x}^j)=c_{j-1}\m{x}^{j-1}$. Three explicit examples of hypermonogenic solutions were computed by choosing the initial function to be a Gau\ss ian function, a Clifford-Bessel function of biaxial type and the Kummer's function respectively.

\section{Cauchy's integral formula for Hypermonogenic solutions in the unit ball}
In this section we derive a Cauchy's formula for hypermonogenic solutions using the Cauchy's formula of monogenic functions. We restrict our study only to the case where the domain is the unit ball in $\mathbb{R}^{p+q}$. 
We assume that a hypermonogenic solution $f$ is defined in a set, which contains a unit ball. The corresponding unit sphere  $S^{p+q-1}\subset\mathbb{R}^{p+q}$  may be expressed  using polar coordinates as a product 
\begin{equation}\label{Nahui}
S^{p+q-1}=S^{p-1}\times S^{q-1}\times [0,\frac{\pi}{2}],
\end{equation}
and then, in these coordinates, a point $\m{\eta}\in S^{p+q-1}$ may be expressed as
\begin{align}\label{Wasa}
\m{\eta}=\cos(\theta)\m{\omega}+\sin(\theta)\m{\nu}
\end{align}
with $\m{\omega}\in S^{p-1}\subset\mathbb{R}^p$, $\m{\nu}\in S^{q-1}\subset\mathbb{R}^q$ and $\theta$ is the angle between vectors $\m{\eta}$ and $\m{\omega}$ with $0\le \theta\le\frac{\pi}{2}$.  We want to point out, that $\m{y}\perp \m{\omega}$, $\m{x}\perp \m{\nu}$ and $\m{\omega}\perp \m{\nu}$. If $|\m{x}+\m{y}|<1$ then the Cauchy formula for a   function $f$ over a unit ball is now
\[
f(\m{x},\m{y})=\frac{1}{\lambda_{p+q-1}}\int_{S^{p+q-1}}\frac{\m{x}+\m{y}-\m{\eta}}{|\m{x}+\m{y}-\m{\eta}|^{p+q}} d\sigma(\m{\eta})f(\m{\eta}),
\]
where $d\sigma(\m{\eta})=\m{\eta}dS(\m{\eta})$. Pulling back the scalar volume element $dS(\m{\eta})$ from the sphere $S^{p+q-1}$ into the preceding decomposition (\ref{Nahui}) via  polar coordinates (\ref{Wasa}), we obtain
\[
dS(\m{\eta})=\cos^{p-1}(\theta)dS(\m{\omega})\sin^{q-1}(\theta)dS(\m{\nu})d\theta
\]
Writing the Cauchy integral using the preceding coordinates, we have
\begin{align*}
&f(\m{x},\m{y})=\\
&\frac{1}{\lambda_{p+q-1}}\int_{0}^{\frac{\pi}{2}}\int_{S^{q-1}}\int_{S^{p-1}}\frac{\m{x}+\m{y}-\cos(\theta)\m{\omega}-\sin(\theta)\m{\nu}}{|\m{x}+\m{y}-\cos(\theta)\m{\omega}-\sin(\theta)\m{\nu}|^{p+q}} (\cos(\theta)\m{\omega}+\sin(\theta)\m{\nu})f(\m{\eta})dS(\m{\omega})d\Gamma(\theta,\m{\nu})
\end{align*}
where $d\Gamma(\theta,\m{\nu})=\cos^{p-1}(\theta)\sin^{q-1}(\theta)dS(\m{\nu})d\theta$. \\
\\
In the preceding case, a hypermonogenic solution 
\[
f(\m{x},\m{y})=A(|\m{x}|,\m{y})+\frac{\m{x}}{|\m{x}|}B(|\m{x}|,\m{y}),
\]
takes the form
\[
f(\m{\eta})=A(\cos(\theta),\sin(\theta)\m{\nu})+ \m{\omega}B(\cos(\theta),\sin(\theta)\m{\nu})
\]
on the sphere. 
We see that funtions $A$ and $B$ do not depend on the variable $\m{\omega}$ and we may simplify the innermost integral. We compute the following two lemmas.

\begin{lem}We may write
\begin{align*}
&\int_{S^{p-1}}\frac{\m{x}+\m{y}-\cos(\theta)\m{\omega}-\sin(\theta)\m{\nu}}{|\m{x}+\m{y}-\cos(\theta)\m{\omega}-\sin(\theta)\m{\nu}|^{p+q}} (\cos(\theta)\m{\omega}+\sin(\theta)\m{\nu})f(\m{\eta})dS(\m{\omega})\\
&=\big(A+(\m{x}+\m{y})(\sin(\theta)\m{\nu}A-\cos(\theta)B)\big)I(\m{x}+\m{y};\theta,\m{\nu})
\end{align*}
where
\[
I(\m{x}+\m{y};\theta,\m{\nu})=\int_{S^{p-1}}\frac{dS(\m{\omega})}{|\m{x}+\m{y}-\cos(\theta)\m{\omega}-\sin(\theta)\m{\nu}|^{p+q}} 
\]
\end{lem}
Proof. Let us first compute
\begin{align*}
&(\m{x}+\m{y}-\cos(\theta)\m{\omega}-\sin(\theta)\m{\nu})(\cos(\theta)\m{\omega}+\sin(\theta)\m{\nu})\\
&=(\m{x}+\m{y}) \big(\cos(\theta)\m{\omega}+\sin(\theta)\m{\nu}\big)+1
\end{align*}
and then
\begin{align*}
&(\m{x}+\m{y}-\cos(\theta)\m{\omega}-\sin(\theta)\m{\nu})(\cos(\theta)\m{\omega}+\sin(\theta)\m{\nu})f(\m{\eta})\\
&=\big((\m{x}+\m{y}) (\cos(\theta)\m{\omega}+\sin(\theta)\m{\nu})+1\big)A\\
&+\big((\m{x}+\m{y}) (-\cos(\theta)+\sin(\theta)\m{\nu}\m{\omega})+\m{\omega}\big)B\\
&=A+(\m{x}+\m{y})(\sin(\theta)\m{\nu}A-\cos(\theta)B)\\
&+(\m{x}+\m{y}) \big(\cos(\theta)\m{\omega}A+ \sin(\theta)\m{\nu}\m{\omega}B\big) +\m{\omega}B.
\end{align*}
Then we obtain
\begin{align*}
&\int_{S^{p-1}}\frac{\m{x}+\m{y}-\cos(\theta)\m{\omega}-\sin(\theta)\m{\nu}}{|\m{x}+\m{y}-\cos(\theta)\m{\omega}-\sin(\theta)\m{\nu}|^{p+q}} (\cos(\theta)\m{\omega}+\sin(\theta)\m{\nu})f(\m{\eta})dS(\m{\omega})\\
&=\int_{S^{p-1}}\frac{A+(\m{x}+\m{y})(\sin(\theta)\m{\nu}A-\cos(\theta)B)}{|\m{x}+\m{y}-\cos(\theta)\m{\omega}-\sin(\theta)\m{\nu}|^{p+q}} )dS(\m{\omega})\\
&+\underbrace{\int_{S^{p-1}}\frac{(\m{x}+\m{y}) \big(\cos(\theta)\m{\omega}A+ \sin(\theta)\m{\nu}\m{\omega}B\big) +\m{\omega}B}{|\m{x}+\m{y}-\cos(\theta)\m{\omega}-\sin(\theta)\m{\nu}|^{p+q}} dS(\m{\omega})}_{=0},
\end{align*}
where the last integral is zero, because the integrand is an odd function with respect to $\m{\omega}$ (see \cite{F}). This proves the lemma.$\square$

\begin{lem} The function $I(\m{x}+\m{y};\theta,\m{\nu})$ in the preceding lemma may be written in the form
\[
I(|\m{x}|,\m{y};\theta,\m{\nu})
= \frac{\lambda_{p-1}}{\tau^{\frac{p+q}{2}}}\frac{\Gamma(p-1)}{2^{p-2}\Gamma^2(\frac{p-1}{2})}\Big(\frac{\tau}{\tau+2r\cos(\theta)}\Big)^{\frac{p+q}{2}}
{}_2F_1\Big(\frac{p+q}{2},\frac{p-1}{2},p-1; \frac{4r\cos(\theta)}{\tau+2r\cos(\theta)}\Big)
\]
where $\tau=|\m{x}|^2+\cos^2(\theta)+|\m{y}-\sin(\theta)\m{\nu}|^2$, $r=|\m{x}|$ and ${}_2F_1$ is a hypergeometric function.
\end{lem}
Proof. Since $\m{x}-\cos(\theta)\m{\omega}\perp \m{y}-\sin(\theta)\m{\nu}$,
we have
\begin{align*}
|\m{x}+\m{y}-\cos(\theta)\m{\omega}-\sin(\theta)\m{\nu}|^2&=|\m{x}-\cos(\theta)\m{\omega}|^2+|\m{y}-\sin(\theta)\m{\nu}|^2\\
&=\tau-2\cos(\theta)\langle \m{x},\m{\omega}\rangle,
\end{align*}
where $\tau=|\m{x}|^2+\cos^2(\theta)+|\m{y}-\sin(\theta)\m{\nu}|^2$. We denote $r=|\m{x}|$ and $\m{\xi}=\frac{\m{x}}{r}\in S^{p-1}$.
Then we obtain, using Funk-Hecke theorem, that
\begin{align*}
I(\m{x}+\m{y};\theta,\m{\nu})&=\int_{S^{p-1}}\frac{dS(\m{\omega})}{\left(\tau-2r\cos(\theta)\langle \m{\xi},\m{\omega}\rangle\right)^{\frac{p+q}{2}}}\\
&= \lambda_{p-1} \int_{-1}^1 \frac{(1-t^2)^{\frac{p-3}{2}}}{(\tau-2r\cos(\theta)t)^{\frac{p+q}{2}}}dt\\
&= \frac{\lambda_{p-1}}{\tau^{\frac{p+q}{2}}} \int_{-1}^1 \frac{(1-t^2)^{\frac{p-3}{2}}}{(1-\frac{2r\cos(\theta)}{\tau}t)^{\frac{p+q}{2}}}dt
\end{align*}
Recall the Euler's formula for  hypergeometric functions (see e.g. \cite{A})
\[
{}_2F_1(a,b,2b; z)=\frac{\Gamma(2b)}{\Gamma^2(b)}\int_0^1s^{b-1}(1-s)^{b-1}(1-zs)^{-a}ds,
\]
which converge for $|z|<1$ and $b>0$. Changing new variables $t= 2s-1$, we obtain the integral
\begin{align*}
{}_2F_1(a,b,2b; z)&=\frac{\Gamma(2b)}{2\Gamma^2(b)}\int_{-1}^1(\frac{1+t}{2})^{b-1}(\frac{1-t}{2})^{b-1}(1-z\frac{t+1}{2})^{-a}dt\\
&=\frac{\Gamma(2b)}{2^{2b-1}(1-\frac{z}{2})^a\Gamma^2(b)}\int_{-1}^1(1-t^2)^{b-1}(1-\frac{z}{2(1-\frac{z}{2})}t)^{-a}dt.
\end{align*}
If we take $a=\frac{p+q}{2}$, $b-1=\frac{p-3}{2}$ and $\frac{z}{2(1-\frac{z}{2})}=\frac{2r\cos(\theta)}{\tau}$ we obtain
\[b=\frac{p-1}{2}, \hspace{1cm} z=\frac{4r\cos(\theta)}{\tau+2r\cos(\theta)}.\]
%\begin{align*}
%b-1=\frac{p-3}{2}\ \Leftrightarrow\ b=\frac{p-1}{2},\\
%\frac{z}{2(1-\frac{z}{2})}=\frac{2r\cos(\theta)}{\tau}\ \Leftrightarrow\ z=4\frac{r\cos(\theta)}{\tau+r\cos(\theta)}
%\end{align*}
Hence,
\begin{multline*}
  \int_{-1}^1(1-t^2)^{\frac{p-3}{2}}(1-\frac{2r\cos(\theta)}{\tau}t)^{-\frac{p+q}{2}}dt \\ = \frac{2^{p-2}\Gamma^2(\frac{p-1}{2})}{\Gamma(p-1)}\Big(\frac{\tau}{\tau+2r\cos(\theta)}\Big)^{\frac{p+q}{2}}
{}_2F_1\Big(\frac{p+q}{2},\frac{p-1}{2},p-1; \frac{4r\cos(\theta)}{\tau+2r\cos(\theta)}\Big) 
\end{multline*}
%\begin{align*}
%&\int_{-1}^1(1-t^2)^{\frac{p-3}{2}}(1-\frac{2r\cos(\theta)}{\tau}t)^{-\frac{p+q}{2}}dt\\
%&=\frac{2^{p-2}\Gamma^2(\frac{p-1}{2})}{\Gamma(p-1)}\Big(\frac{\tau}{\tau+2r\cos(\theta)}\Big)^{\frac{p+q}{2}}
%{}_2F_1\Big(\frac{p+q}{2},\frac{p-1}{2},p-1; 4\frac{r\cos(\theta)}{\tau+r\cos(\theta)}\Big)
%\end{align*}
Then we have
\begin{align*}
I(\m{x}+\m{y};\theta,\m{\nu})
&= \frac{\lambda_{p-1}}{\tau^{\frac{p+q}{2}}}\frac{2^{p-2}\Gamma^2(\frac{p-1}{2})}{\Gamma(p-1)}\Big(\frac{\tau}{\tau+2r\cos(\theta)}\Big)^{\frac{p+q}{2}}
{}_2F_1\Big(\frac{p+q}{2},\frac{p-1}{2},p-1; \frac{4r\cos(\theta)}{\tau+2r\cos(\theta)}\Big),
\end{align*}
which depends only on $|\m{x}|$ and $\m{y}$.$\square$\\
\\
Using the preceding lemmas, we obtain the following Cauchy's formula
\begin{align*}
&f(\m{x},\m{y})=\\
&\frac{1}{\lambda_{p+q-1}}\int_{0}^{\frac{\pi}{2}}\int_{S^{q-1}}
\big(A+(\m{x}+\m{y})(\sin(\theta)\m{\nu}A-\cos(\theta)B)\big)I(|\m{x}|,\m{y};\theta,\m{\nu})
d\Gamma(\theta,\m{\nu})\\
&=\frac{1}{\lambda_{p+q-1}}\int_{0}^{\frac{\pi}{2}}\int_{S^{q-1}}I(|\m{x}|,\m{y};\theta,\m{\nu})
(1+(\m{x}+\m{y})\sin(\theta)\m{\nu})A(\cos(\theta),\sin(\theta)\m{\nu})
d\Gamma(\theta,\m{\nu})\\
&-\frac{1}{\lambda_{p+q-1}}\int_{0}^{\frac{\pi}{2}}\int_{S^{q-1}}I(|\m{x}|,\m{y};\theta,\m{\nu})
\cos(\theta)B(\cos(\theta),\sin(\theta)\m{\nu})
d\Gamma(\theta,\m{\nu})
\end{align*}
Since $f(\m{x},\m{y})=A(|\m{x}|,\m{y})+\frac{\m{x}}{|\m{x}|}B(|\m{x}|,\m{y})$ we obtain the Cauchy's formulas for $A$ and $B$ parts of the function separately.

\begin{thm} In the unit ball
\[
A(|\m{x}|,\m{y})=\frac{1}{\lambda_{p+q-1}}\int_{0}^{\frac{\pi}{2}}\int_{S^{q-1}}I(|\m{x}|,\m{y};\theta,\m{\nu})
\big(A+\m{y}(\sin(\theta)\m{\nu}A-\cos(\theta)B)\big)
d\Gamma(\theta,\m{\nu})
\]
and
\[
B(|\m{x}|,\m{y})=\frac{|\m{x}|}{\lambda_{p+q-1}}\int_{0}^{\frac{\pi}{2}}\int_{S^{q-1}}
I(|\m{x}|,\m{y};\theta,\m{\nu})\sin(\theta)\m{\nu}A\,d\Gamma(\theta,\m{\nu})
\]
where $d\Gamma(\theta,\m{\nu})=\cos^{p-1}(\theta)\sin^{q-1}(\theta)dS(\m{\nu})d\theta$. 
\end{thm}

\section{Hypermonogenic Plane Wave Solutions}
In this section our aim is to give a definition for hypermonogenic plane wave solutions. In general, plane wave solutions in Clifford analysis are solutions of the Dirac operator depending only on inner products of the variable with the unit normal vector of  the wave's direction of propagation. Monogenic plane waves are allways complex valued functions of the form
\[
f(\langle \m{x},\m{s}\rangle)\m{s},
\]
where $\m{s}^2=\m{0}$ and $f(z)$ is a holomorphic function.
To find a right definition for hypermonogenic plane waves,  we consider, as an example,  a biaxial monogenic plane wave
\[
h(\m{x},\m{y})=\langle \m{x}+\m{y},\m{\tau}+\m{\sigma}\rangle^k(\m{\tau}+\m{\sigma})
\] 
where $\m{\tau}\in\mathbb{C}^p$ and $\m{\sigma}\in\mathbb{C}^q$. This function is monogenic if $\m{\tau}+\m{\sigma}$ is a null vectors, i.e., 
\[
(\m{\tau}+\m{\sigma})^2=\m{\tau}^2+\m{\sigma}^2=0.
\]
We may take $\m{\tau}=\m{t}$ and $\m{\sigma}=i\m{s}$, where $\m{t}\in S^{p-1}$ and $\m{s}\in S^{q-1}$. Then we have
\[
h(\m{x},\m{y})=(\langle \m{x},\m{t}\rangle+i\langle\m{y},\m{s}\rangle)^k(\m{t}+i\m{s})
\]
Since hypermonogenic solutions are radial with respect to $\m{x}$, we may radialize the function $h$ by integrating it over $\m{t}\in S^{p-1}$, i.e., we define the function
\[
g(\m{x},\m{y})=\int_{S^{p-1}} (\langle \m{x},\m{t}\rangle+i\langle\m{y},\m{s}\rangle)^k(\m{t}+i\m{s})dS(\m{t}).
\]
Now we want to see the proper form of this integral for a definition. We start from the following lemma.

\begin{lem}
There exist functions $A$ and $B$ such that
\[
g(\m{x},\m{y})=A+iB\m{s}.
\]
\end{lem}
Proof. We compute
\begin{align*}
g(\m{x},\m{y})&=\int_{S^{p-1}} (\langle \m{x},\m{t}\rangle+i\langle\m{y},\m{s}\rangle)^k(\m{t}+i\m{s})dS(\m{t})\\
&=\int_{S^{p-1}} (\langle \m{x},\m{t}\rangle+i\langle\m{y},\m{s}\rangle)^k\m{t}dS(\m{t})\\
&+i\int_{S^{p-1}} (\langle \m{x},\m{t}\rangle+i\langle\m{y},\m{s}\rangle)^k dS(\m{t})\m{s}.
\end{align*}
\begin{flushright}
$\square$
\end{flushright}
We obtain the following explicit expressions for these functions.

\begin{prop} The function $A$ is of the form
\[
A=\lambda_{p-1}\sum_{j=0}^{[\frac{k-1}{2}]}  \m{x}^{2j+1}  a_j(\langle\m{y},\m{s}\rangle)
\]
where
\[
a_j(\langle\m{y},\m{s}\rangle)=(-1)^j  {k \choose 2j+1}\frac{\Gamma(\frac{p-1}{2})\Gamma(j+\frac{3}{2})}{\Gamma(\frac{p}{2}+j+1)}(i\langle\m{y},\m{s}\rangle)^{k-2j-1}.
\]
\end{prop}
Proof. Using binomial theorem we have
\begin{align*}
A&=\int_{S^{p-1}} (\langle \m{x},\m{t}\rangle+i\langle\m{y},\m{s}\rangle)^k\m{t}dS(\m{t})\\
&=\int_{S^{p-1}} \sum_{j=0}^k {k \choose j}  \langle \m{x},\m{t}\rangle^j(i\langle\m{y},\m{s}\rangle)^{k-j}\m{t}dS(\m{t})\\
&=\sum_{j=0}^k {k \choose j} (i\langle\m{y},\m{s}\rangle)^{k-j} \int_{S^{p-1}}   \langle \m{x},\m{t}\rangle^j\m{t}dS(\m{t})
\end{align*}
We recall that for odd integrands the above integral is zero (see \cite{F}), i.e., 
\[
\int_{S^{p-1}}   \langle \m{x},\m{t}\rangle^{2j}\m{t}dS(\m{t})=0.
\]
For even integrands we use the Funk-Hecke formula with $H_1(\m{x})=\m{x}$, $\psi(u)=u^{2j+1}$ and $\frac{1!}{(p-2)_1} C^{\frac{p}{2}-1}_1(u)=u$. Hence,
\begin{align*}
\int_{S^{p-1}}   \langle \m{x},\m{t}\rangle^{2j+1}\m{t}dS(\m{t})&=|\m{x}|^{2j+1}\int_{S^{p-1}}   \langle \frac{\m{x}}{|\m{x}|},\m{t}\rangle^{2j+1}\m{t}dS(\m{t})\\
&=|\m{x}|^{2\ell+1} \frac{\m{x}}{|\m{x}|} \lambda_{p-1}\int_{-1}^{1} u^{2j+1} u(1-u^2)^{\frac{p-3}{2}}du.
\end{align*}
 Then as a routine integration we obtain
\[
\int_{-1}^{1} u^{2j+1} u(1-u^2)^{\frac{p-3}{2}}du=\frac{\Gamma(\frac{p-1}{2})\Gamma(j+\frac{3}{2})}{\Gamma(\frac{p}{2}+j+1)},
\]
that is
\[
\int_{S^{p-1}}   \langle \m{x},\m{t}\rangle^{2j+1}\m{t}dS(\m{t})=(-1)^j\m{x}^{2j+1}  \lambda_{p-1}\frac{\Gamma(\frac{p-1}{2})\Gamma(j+\frac{3}{2})}{\Gamma(\frac{p}{2}+j+1)}.
\]
Then
\begin{align*}
A&=\sum_{j=0}^{[\frac{k-1}{2}]} {k \choose 2j+1}\int_{S^{p-1}}   \langle \m{x},\m{t}\rangle^{2j+1}\m{t}dS(\m{t})(i\langle\m{y},\m{s}\rangle)^{k-2j-1}\\
&=\sum_{j=0}^{[\frac{k-1}{2}]}  \m{x}^{2j+1}  {k \choose 2j+1}(-1)^j  \lambda_{p-1}\frac{\Gamma(\frac{p-1}{2})\Gamma(j+\frac{3}{2})}{\Gamma(\frac{p}{2}+j+1)}(i\langle\m{y},\m{s}\rangle)^{k-2j-1}.
\end{align*}
\begin{flushright}
$\square$
\end{flushright}

\begin{prop} The function $B$ is of the form
\[
B=\lambda_{p-1}\sum_{j=0}^{[\frac{k}{2}]} \m{x}^{2j}b_j(\langle\m{y},\m{s}\rangle)
\]
where
\[
b_j(\langle\m{y},\m{s}\rangle)=(-1)^j{k \choose 2j}\frac{\Gamma(\frac{p-1}{2})\Gamma(j+\frac{1}{2})}{\Gamma(\frac{p}{2}+j)}(i\langle\m{y},\m{s}\rangle)^{k-2j}.
\]
\end{prop}
Proof. Similarly as in the preceding proposition, we compute
\begin{align*}
B&=\int_{S^{p-1}} (\langle \m{x},\m{t}\rangle+i\langle\m{y},\m{s}\rangle)^k dS(\m{t})\\
&=\sum_{j=0}^k {k \choose j}\int_{S^{p-1}}   \langle \m{x},\m{t}\rangle^jdS(\m{t})(i\langle\m{y},\m{s}\rangle)^{k-j}
\end{align*}
and since $\m{t}\mapsto\langle \m{x},\m{t}\rangle^{2j+1}$ is an odd function we obtain
\[
\int_{S^{p-1}}   \langle \m{x},\m{t}\rangle^{2j+1}dS(\m{t})=0.
\]
Using Funk-Hecke formula, where $
H_0(\m{x})=1,\ \frac{0!}{(p-2)_0}C^{\frac{p}{2}-1}_0(u)=1
$, we compute
\begin{align*}
&\int_{S^{p-1}}   \langle \m{x},\m{t}\rangle^{2j}dS(\m{t})\\
&=|\m{x}|^{2j}\int_{S^{p-1}}   \langle \frac{\m{x}}{|\m{x}|},\m{t}\rangle^{2j}dS(\m{t})\\
&=|\m{x}|^{2j}\lambda_{p-1}\int_{-1}^{1} u^{2j} (1-u^2)^{\frac{p-3}{2}}du\\
&=|\m{x}|^{2j}\lambda_{p-1}\frac{\Gamma(\frac{p-1}{2})\Gamma(j+\frac{1}{2})}{\Gamma(\frac{p}{2}+j)}.
\end{align*}
Substituting these to the original formula, we have
\begin{align*}
B&=\int_{S^{p-1}} (\langle \m{x},\m{t}\rangle+i\langle\m{y},\m{s}\rangle)^k dS(\m{t})\\
&=\sum_{j=0}^{[\frac{k}{2}]} {k \choose 2j}|\m{x}|^{2j}\lambda_{p-1}\frac{\Gamma(\frac{p-1}{2})\Gamma(j+\frac{1}{2})}{\Gamma(\frac{p}{2}+j)}(i\langle\m{y},\m{s}\rangle)^{k-2j}\\
&=\sum_{j=0}^{[\frac{k}{2}]} \m{x}^{2j}(-1)^j\lambda_{p-1}{k \choose 2j}\frac{\Gamma(\frac{p-1}{2})\Gamma(j+\frac{1}{2})}{\Gamma(\frac{p}{2}+j)}(i\langle\m{y},\m{s}\rangle)^{k-2j}
\end{align*}
$\square$\\
\\
We see that we may write the monogenic function $g$ of the form
\[
g(\m{x},\m{y})=\lambda_{p-1} \sum_{j=0}^{[\frac{k-1}{2}]}  \m{x}^{2j+1}  a_j(\langle\m{y},\m{s}\rangle)+\lambda_{p-1}\sum_{j=0}^{[\frac{k}{2}]} \m{x}^{2j}ib_j(\langle\m{y},\m{s}\rangle)\m{s}
\]
This motivates us to define hypermonogenic plane waves as monogenic functions of the form 
\begin{align}\label{HPW}
F(\m{x},\langle\m{y},\m{s}\rangle)=\sum_{j=0}^\infty \m{x}^j (C_j(\langle\m{y},\m{s}\rangle)+\m{s}D_j(\langle\m{y},\m{s}\rangle)),
\end{align}
where functions $C_j$ and $D_j$ are complex valued differentiable functions of one complex variable. This definition is a good starting point for the construction of the hypermonogenic plane waves and we will exploit it later. Developing this series expansion in terms of powers of $|\underline{x}|$, we see that functions are always of the following form.

\begin{defn}[Hypermonogenic plane waves (HPW)]
They are functions of the form
\[
F(\m{x},\langle\m{y},\m{s}\rangle)=A(|\m{x}|,\langle\m{y},\m{s}\rangle)+B(|\m{x}|,\langle\m{y},\m{s}\rangle)\frac{\m{x}}{|\m{x}|}+C(|\m{x}|,\langle\m{y},\m{s}\rangle)\m{s}+D(|\m{x}|,\langle\m{y},\m{s}\rangle)\frac{\m{x}}{|\m{x}|}\m{s}
\]
where the coefficient functions $A$, $B$, $C$ and $D$ are differentiable complex valued functions, $s\in S^{q-1}$ and   $(\partial_{\m{x}}+\partial_{\m{y}})F=0$.
\end{defn}
It is obvious, that the hypermonogenic plane wave solutions are hypermonogenic solutions.\\
\\
Next we derive a system for coeffecient functions of a HPW of the form (\ref{HPW}) and for that we recall
\[
\partial_{\m{x}} \m{x}^j=\begin{cases}
-j\m{x}^{j-1}, & \text{ for $j$ even},\\ 
-(j+p-1)\m{x}^{j-1}, & \text{ for $j$ odd},\\ 
\end{cases}
\]
where we may briefly denote
\[
\partial_{\m{x}} \m{x}^j=\beta_j\m{x}^{j-1}.
\]

\begin{prop}\label{P5}
A function  $F(\m{x},\langle\m{y},\m{s}\rangle)$ of the form (\ref{HPW}) is a hypermonogenic plane wave if and only if
\begin{align*}
\beta_{j+1}C_{j+1}(t) -(-1)^j D'_j(t)=0,\\
\beta_{j+1}D_{j+1}(t) +(-1)^j C'_j(t)=0,
\end{align*}
where  $t=\langle\m{y},\m{s}\rangle$ and it is determined by its initial functions $C_0(t)$ and $D_0(t)$.
\end{prop}
Proof. First we have
\[
\partial_{\m{y}}C_j(\langle\m{y},\m{s}\rangle)=\m{s}C'_j(\langle\m{y},\m{s}\rangle).
\]
Since $\partial_{\m{y}}\m{x}^j=(-1)^{j}\m{x}^j\partial_{\m{y}}$, we obtain
\begin{align*}
(\partial_{\m{x}}+\partial_{\m{y}})F(\m{x},\langle\m{y},\m{s}\rangle)
&=\sum_{j=1}^\infty \beta_j\m{x}^{j-1}(C_j(\langle\m{y},\m{s}\rangle)+\m{s}D_j(\langle\m{y},\m{s}\rangle))\\
&+\sum_{j=0}^\infty \m{x}^j (-1)^j(\m{s}C'_j(\langle\m{y},\m{s}\rangle)+D'_j(\langle\m{y},\m{s}\rangle)\m{s}^2)\\
&=\sum_{j=0}^\infty \m{x}^{j}(\beta_{j+1}C_{j+1}(\langle\m{y},\m{s}\rangle)-(-1)^jD'_j(\langle\m{y},\m{s}\rangle))\\
&+\sum_{j=0}^\infty \m{x}^j (\beta_{j+1}D_{j+1}(\langle\m{y},\m{s}\rangle)+ (-1)^jC'_j(\langle\m{y},\m{s}\rangle))\m{s}\\
\end{align*}
completing the proof.$\square$

\section{Two methods to find hypermonogenic plane waves}
In this section, we describe two methods that show how to construct hypermonogenic plane wave solutions. \textbf{The first method} uses series of the  form (\ref{HPW}) and we illustrate it with the following example.

\begin{prop}\label{JuriS}
If we look for a HPW of the form
\[
F(\m{x},\langle\m{y},\m{s}\rangle)=\sum_{j=0}^\infty \m{x}^j (c_j +\m{s}d_j)e^{\langle\m{y},\m{s}\rangle}
\]
with the initial condition $c_0=1$ and $d_0=0$ we obtain the solution
\[
F(\m{x},\langle\m{y},\m{s}\rangle)=\frac{2^{\frac{p}{2}-1}\Gamma(\frac{p}{2})}{|\m{x}|^{\frac{p}{2}-1}}\Big(J_{\frac{p}{2}-1}(|\m{x}|)  +J_{\frac{p}{2}}(|\m{x}|) \frac{\m{x}}{|\m{x}|}\m{s}\Big)e^{\langle\m{y},\m{s}\rangle}
\]
where $J_{\nu}$ are Bessel functions.
\end{prop}
Proof. This particular case can be written in the form (\ref{HPW}) by taking
\[C_j(t)=c_je^{t}, \;\;D_j(t)=d_je^{t}, \;\;\;\;\mbox{ which implies }\;\;\;\; C'_j(t)=c_je^{t}, \;\;D'_j(t)=d_je^{t}.\]
%\begin{align*}
%C_j(t)=c_je^{t}, \ C'_j(t)=c_je^{t},\\
%D_j(t)=d_je^{t},\ D'_j(t)=d_je^{t}.
%\end{align*}
The system  given in Proposition \ref{P5} becomes now purely algebraic
\[ \begin{cases} \beta_{j+1}c_{j+1} -(-1)^j d_j=0,\\ \beta_{j+1}d_{j+1} +(-1)^j c_j=0, \end{cases} \;\;\;\;\mbox{ or equivalently }\;\;\;\; \begin{cases}
c_{j+1} =\frac{(-1)^j}{\beta_{j+1}} d_j,\\d_{j+1} =-\frac{(-1)^j}{\beta_{j+1}} c_j. \end{cases}\]

%\begin{align*}
%\beta_{j+1}c_{j+1} -(-1)^j d_j=0,\\
%\beta_{j+1}d_{j+1} +(-1)^j c_j=0,
%\end{align*}
%or
%\begin{align*}
%c_{j+1} =\frac{(-1)^j}{\beta_{j+1}} d_j,\\
%d_{j+1} =-\frac{(-1)^j}{\beta_{j+1}} c_j.
%\end{align*}
Since $d_0=0$ we compute that
\[
0=d_0=c_1=d_2=c_3=...
\]
i.e. even $d_j$'s and odd $c_j$'s are zero. Starting from $c_0=1$ we get
\begin{align*}
d_{1} &=-\frac{(-1)^0}{\beta_{1}} c_0=\frac{1}{p}, \\
c_{2} &=\frac{(-1)^1}{\beta_{2}} d_1= \frac{1}{2\cdot p},\\
d_{3} &=-\frac{(-1)^2}{\beta_{3}} c_2=\frac{1}{2\cdot p (p+2)},\\
c_{4} &=\frac{(-1)^3}{\beta_{4}} d_3=\frac{1}{2\cdot 4\cdot p (p+2)},\\
d_{5} &=-\frac{(-1)^4}{\beta_{5}} c_4=\frac{1}{2\cdot 4\cdot p (p+2)(p+4)},\\
c_{6}& =\frac{(-1)^5}{\beta_{6}} d_5=\frac{1}{2\cdot 4\cdot 6\cdot p (p+2)(p+4)}.
\end{align*}
Then we easily infer that when $j$ is even
\[
c_j=\frac{1}{2\cdot 4\cdots j\cdot p(p+2)\cdots (p+j-2)}
\]
and when $j$ is odd
\[
d_j=\frac{1}{2\cdot 4\cdots (j-1)\cdot p\cdot (p+2)\cdots (p+j-1)}.
\]
\underline{If $j$ is even}, we observe that $2^\frac{j}{2} (j/2)!=2\cdot 4\cdots j$ and we may write
\[
c_j=\frac{1}{2^\frac{j}{2} (j/2)!\, p(2+p)\cdots (j-2+p)}
\]
where in addition $(j/2)!=\Gamma(j/2+1)$.
In the product $p(2+p)\cdots (j-2+p)$, there are $j/2$ elements. Then we have
\begin{align*}
&p(2+p)(4+p)\cdots (j-2+p)\\
&=2^{j/2}(0+\frac{p}{2})(1+\frac{p}{2})(2+\frac{p}{2})\cdots (\frac{j}{2}-1+\frac{p}{2})\\
&=2^{j/2}\prod_{k=1}^{j/2} \Big( k-1 +\frac{p}{2}\Big)=2^{j/2}\frac{\Gamma(j/2+p/2)}{\Gamma(p/2)}
\end{align*}
where in the last line we use a classical identity of the Gamma functions. Hence we obtain
\begin{align*}
c_j&=  \frac{\Gamma(p/2)}{2^j \Gamma(j/2+1)\Gamma(j/2+p/2)}.
\end{align*}
\underline{When $j$ is odd}, we observe that $2^{\frac{j-1}{2}} (\frac{j-1}{2})!=2\cdot 4\cdots (j-1)$. Then we have
\[
d_j=\frac{1}{2^{\frac{j-1}{2}} (\frac{j-1}{2})!\cdot p\cdot (2+p)\cdots (j-1+p)}.
\]
where $(\frac{j-1}{2})!=\Gamma(\frac{j+1}{2})$. 
In the product  $p(2+p)(4+p)\cdots (j-1+p)$ there are $\frac{j+1}{2}$ elements.
 Then we have
\begin{align*}
&p(2+p)(4+p)\cdots (j-1+p)\\
&=2^{\frac{j+1}{2}}(0+\frac{p}{2})(1+\frac{p}{2})(2+\frac{p}{2})\cdots (\frac{j-1}{2}+\frac{p}{2})\\
&=2^{\frac{j+1}{2}}\prod_{k=1}^{\frac{j+1}{2}}\Big(k-1+\frac{p}{2}\Big)=2^{\frac{j+1}{2}}\frac{\Gamma((j+1)/2+p/2)}{\Gamma(p/2)}.
\end{align*}
As a consequence,
\[
d_j=\frac{\Gamma(p/2)}{2^{j} \Gamma(\frac{j+1}{2})\Gamma((j+1)/2+p/2)}
\]
This way we obtain the solution
\begin{align*}
F(\m{x},\langle\m{y},\m{s}\rangle)
&=\Gamma(\frac{p}{2})\sum_{j=0}^\infty  \left( \frac{\m{x}^{2j}}{2^{2j} \Gamma(j+1)\Gamma(j+\frac{p}{2})} 
 +  \frac{\m{x}^{2j+1}}{2^{2j+1} \Gamma(j+1)\Gamma(j+1+\frac{p}{2})}\m{s}\right)e^{\langle\m{y},\m{s}\rangle}\\
&=\Gamma(\frac{p}{2})\sum_{j=0}^\infty  \left( \frac{(-1)^j |\m{x}|^{2j}}{2^{2j} \Gamma(j+1)\Gamma(j+\frac{p}{2})} +\frac{(-1)^j |\m{x}|^{2j}}{2^{2j+1} \Gamma(j+1)\Gamma(j+1+\frac{p}{2})}\m{x}\m{s}\right)e^{\langle\m{y},\m{s}\rangle}
\end{align*}
Since $\Gamma(j+1)=j!$ we obtain
\begin{align}
F(\m{x},\langle\m{y},\m{s}\rangle)
&=\Gamma(\frac{p}{2})\left(\sum_{j=0}^\infty   \frac{(-1)^j |\m{x}|^{2j}}{2^{2j} j!\Gamma(\frac{p}{2}+j)}
 +\frac{1}{2}\sum_{j=0}^\infty \frac{(-1)^j |\m{x}|^{2j}}{2^{2j} j!\Gamma(\frac{p}{2}+j+1)}\m{x}\m{s}\right)e^{\langle\m{y},\m{s}\rangle}.\label{UNTERMENSCH}
\end{align}
The sums are now precisely the same as in the definition of a Bessel function
\begin{align}\label{MashinaVremeni}
J_{\nu}(z)=\frac{z^\nu}{2^\nu}\sum_{j=0}^\infty\frac{(-1)^j z^{2j}}{j! 2^{2j}\Gamma(\nu+j+1)}. 
\end{align}
and we obtain the corresponding solution
\begin{align*}
F(\m{x},\langle\m{y},\m{s}\rangle)
&=\Gamma(\frac{p}{2})\left(\frac{2^{\frac{p}{2}-1}}{|\m{x}|^{\frac{p}{2}-1}}J_{\frac{p}{2}-1}(|\m{x}|)  +\frac{1}{2}\frac{2^{\frac{p}{2}}}{|\m{x}|^{\frac{p}{2}}}J_{\frac{p}{2}}(|\m{x}|) \m{x}\m{s}\right)e^{\langle\m{y},\m{s}\rangle}\\
&=\frac{2^{\frac{p}{2}-1}\Gamma(\frac{p}{2})}{|\m{x}|^{\frac{p}{2}-1}}\left(J_{\frac{p}{2}-1}(|\m{x}|)  +J_{\frac{p}{2}}(|\m{x}|) \frac{\m{x}}{|\m{x}|}\m{s}\right)e^{\langle\m{y},\m{s}\rangle}.
\end{align*}
\begin{flushright}
$\square$
\end{flushright}
The preceding proof is a direct application of the fundamental system given in Proposition \ref{P5}. The same solution may be found using the generalized Cauchy-Kovalevskaya extension. 

\begin{prop}
The plane wave solution found in Proposition \ref{JuriS} may be obtained using the generalized Cauchy-Kovalevskaya extension with the initial function $e^{\langle\m{y},\m{s}\rangle}$.
\end{prop}
Proof. Using the definition of a Bessel function (\ref{MashinaVremeni}) we find the formulas
\begin{align*}
\Big(\frac{|\m{x}|}{2}\sqrt{\Delta_{\m{y}}}\Big)^{-p/2}\Big(\frac{|\m{x}|}{2}\sqrt{\Delta_{\m{y}}}J_{\frac{p}{2}-1}(|\m{x}|\sqrt{\Delta_{\m{y}}})\Big)&=\sum_{j=0}^\infty \frac{(-1)^j |\m{x}|^{2j} \Delta_{\m{y}}^j}{j! 2^{2j}\Gamma(\frac{p}{2}+j)},\\
\Big(\frac{|\m{x}|}{2}\sqrt{\Delta_{\m{y}}}\Big)^{-p/2}J_{\frac{p}{2}}(|\m{x}|\sqrt{\Delta_{\m{y}}})&=\sum_{j=0}^\infty \frac{(-1)^j |\m{x}|^{2j} \Delta_{\m{y}}^j}{j! 2^{2j}\Gamma(\frac{p}{2}+j+1)},
\end{align*}
which we need in the formula given in Theorem \ref{Absolut}.
Since $\partial_{\m{y}}e^{\langle\m{y},\m{s}\rangle}=e^{\langle\m{y},\m{s}\rangle}\m{s}$ we compute
\[
\Delta_{\m{y}}e^{\langle\m{y},\m{s}\rangle}=-\partial_{\m{y}}^2e^{\langle\m{y},\m{s}\rangle}=-\partial_{\m{y}}e^{\langle\m{y},\m{s}\rangle}\m{s}=
-e^{\langle\m{y},\m{s}\rangle}\m{s}^2=e^{\langle\m{y},\m{s}\rangle}, 
\]
and then $\Delta^j_{\m{y}}e^{\langle\m{y},\m{s}\rangle}=e^{\langle\m{y},\m{s}\rangle}$ for all $j=0,1,2,...$.
Using the generalized Cauchy-Kovalevskaya formula we obtain
\begin{align*}
f(\m{x},\m{y})&=\Gamma(\frac{p}{2})\Big(\frac{|\m{x}|}{2}\sqrt{\Delta_{\m{y}}}\Big)^{-p/2}\Big(\frac{|\m{x}|}{2}\sqrt{\Delta_{\m{y}}}J_{\frac{p}{2}-1}(|\m{x}|\sqrt{\Delta_{\m{y}}})+\frac{\m{x}\partial_{\m{y}}}{2}J_{\frac{p}{2}}(|\m{x}|\sqrt{\Delta_{\m{y}}})\Big)e^{\langle\m{y},\m{s}\rangle}\\
&=\Gamma(\frac{p}{2})\Big(\sum_{j=0}^\infty \frac{(-1)^j |\m{x}|^{2j} }{j! 2^{2j}\Gamma(\frac{p}{2}+j)} +\frac{1}{2}\sum_{j=0}^\infty \frac{(-1)^j |\m{x}|^{2j}}{j! 2^{2j}\Gamma(\frac{p}{2}+j+1)}\m{x}\m{s}\Big)e^{\langle\m{y},\m{s}\rangle},
\end{align*}
which coincides with the function obtained in (\ref{UNTERMENSCH}).$\square$\\
\\
\textbf{The second method} to construct hypermonogenic plane wave solutions is based on the Funk-Hecke formula, similarly as in the beginning of this section. We will start from a monogenic plane wave,  let us take for example
\begin{align}\label{Teroesa}
e^{\langle \m{x},\m{t}\rangle +i\langle \m{y},\m{s}\rangle}(\m{t}+i\m{s}).
\end{align}
Since a hypermongenic plane wave must be a radial function with respect to the variable $\m{x}$, we need to integrate over the sphere $\m{t}\in S^{p-1}$ and then we have the function
\[
G(\m{x},\langle \m{y},\m{s}\rangle)=\int_{S^{p-1}}e^{\langle \m{x},\m{t}\rangle +i\langle \m{y},\m{s}\rangle}(\m{t}+i\m{s})dS(\m{t}).
\]
Nowit suffices to compute the integral using the Funk-Hecke formula, which we will do in the next proposition. As a result we obtain a hypermonogenic plane wave solution. 

This example demonstrate the general idea of the method. Making the summary, instead of function (\ref{Teroesa}), we may start from an any monogenic plane wave and radialize it integrating with respect to $\m{t}$. This integral may be computed by using Funk-Hecke formula. The resulting function will be a hypermonogenic plane wave.

\begin{prop}
The preceding function $G(\m{x},\langle \m{y},\m{s}\rangle)$ is a hypermonogenic plane wave and can be written as
\[
G(\m{x},\langle \m{y},\m{s}\rangle)=\frac{\sqrt{\pi}\lambda_{p-1}2^\frac{p-2}{2}\Gamma(\frac{p-1}{2})}{|\m{x}|^{\frac{p-2}{2}}} \Big(iI_\frac{p-2}{2}(|\m{x}|)\m{s}+ \frac{\m{x}}{|\m{x}|} I_{\frac{p}{2}}(|\m{x}|)
\Big)e^{i\langle \m{y},\m{s}\rangle},
\]
where $I$ is a modified Bessel function.
\end{prop}
Proof. We observe that we may write $G$ in the form 
\[
G(\m{x},\langle \m{y},\m{s}\rangle)=A(\m{x},\langle \m{y},\m{s}\rangle)+iB(\m{x},\langle \m{y},\m{s}\rangle)\m{s}
\]
where
\[
A=e^{i\langle \m{y},\m{s}\rangle}\int_{S^{p-1}}e^{\langle \m{x},\m{t}\rangle} \m{t}dS(\m{t})
\]
and
\[
B=e^{i\langle \m{y},\m{s}\rangle}\int_{S^{p-1}}e^{\langle \m{x},\m{t}\rangle} dS(\m{t}).
\]
Let us now compute the integrals. We define $\m{x}=r\m{v}$, where $\m{v}\in S^{p-1}$.  Funk-Hecke formula gives us
\begin{align*}
&\int_{S^{p-1}}e^{r\langle \m{v},\m{t}\rangle} \m{t}dS(\m{t})=\lambda_{p-1}\m{v}\int_{-1}^{1}e^{ru}u(1-u^2)^{\frac{p-3}{2}}du
\end{align*}
and 
\begin{align*}
&\int_{S^{p-1}}e^{\langle \m{x},\m{t}\rangle} dS(\m{t})
=\lambda_{p-1}\int_{-1}^{1}e^{ru}(1-u^2)^{\frac{p-3}{2}}du
\end{align*}
We recall the integral representation to a modified Bessel function ($\text{Re}(\nu)>-\frac{1}{2}$)
\[
I_\nu(z)=\frac{1}{\Gamma(\nu+\frac{1}{2})\sqrt{\pi}} \big(\frac{z}{2}\big)^\nu\int_{-1}^1 e^{-zu}(1-u^2)^{\nu-\frac{1}{2}}du
\]
If we put $\nu=\frac{p-2}{2}$ and $r=-z$, we obtain 
\begin{align*}
&I_\frac{p-2}{2}(-r)=(-1)^\frac{p-2}{2}\frac{1}{\Gamma(\frac{p-1}{2})\sqrt{\pi}} \big(\frac{r}{2}\big)^\frac{p-2}{2}\int_{-1}^1 e^{ru}(1-u^2)^{\frac{p-3}{2}}du\\
\Leftrightarrow\ &    \int_{-1}^1 e^{ru}(1-u^2)^{\frac{p-3}{2}}du=(-1)^\frac{p-2}{2} \Gamma(\frac{p-1}{2})\sqrt{\pi}  \big(\frac{2}{r}\big)^\frac{p-2}{2} I_\frac{p-2}{2}(-r)
\end{align*}
In order to get rid of the minus-sign, we may use the following well known identity
\[
I_\nu(-z)=(-z)^\nu z^{-\nu}I_\nu(z).
\]
what gives us
\[
I_\frac{p-2}{2}(-r)=(-1)^\frac{p-2}{2}I_\frac{p-2}{2}(r)
\]
and then 
\begin{align*}
 \int_{-1}^1 e^{ru}(1-u^2)^{\frac{p-3}{2}}du= \Gamma(\frac{p-1}{2})\sqrt{\pi}  \big(\frac{2}{r}\big)^\frac{p-2}{2} I_\frac{p-2}{2}(r).
\end{align*}
Subsituting $r=|\m{x}|$ we obtain
\[
B=\lambda_{p-1}\Gamma(\frac{p-1}{2})\sqrt{\pi}e^{i\langle \m{y},\m{s}\rangle}\big(\frac{2}{|\m{x}|}\big)^{\frac{p-2}{2}}I_\frac{p-2}{2}(|\m{x}|).
\]
Consider now a function
\[
\varphi(r)=\int_{-1}^{1}e^{ru}(1-u^2)^{\frac{p-3}{2}}du
\]
Since the integrated function and its derivative with respect to $r$ are continuous, we may compute
\[
\varphi'(r)=\int_{-1}^{1}e^{ru}u(1-u^2)^{\frac{p-3}{2}}du.
\]
Let us recall the following differentiation formula for modified Bessel functions
\[
\frac{d }{dz}(z^{-\nu}I_\nu(z))=z^{-\nu}I_{\nu+1}(z).
\]
Using these, we have
\begin{align*}
\int_{-1}^1 e^{ru}u(1-u^2)^{\frac{p-3}{2}}du&= 2^\frac{p-2}{2}\Gamma(\frac{p-1}{2})\sqrt{\pi}  \frac{d}{dr}\Big(r^{-\frac{p-2}{2}} I_\frac{p-2}{2}(r)\Big)\\
&=\Gamma(\frac{p-1}{2})\sqrt{\pi}   \big(\frac{2}{r}\big)^\frac{p-2}{2} I_{\frac{p}{2}}(r).
\end{align*}
Substituting $r=|\m{x}|$ and $\m{v}=\frac{\m{x}}{|\m{x}|}$ we obtain
\begin{align*}
A&=\lambda_{p-1}e^{i\langle \m{y},\m{s}\rangle}\m{v}\int_{-1}^{1}e^{ru}u(1-u^2)^{\frac{p-3}{2}}du\\
&=\lambda_{p-1}e^{i\langle \m{y},\m{s}\rangle}\frac{\m{x}}{|\m{x}|}\Gamma(\frac{p-1}{2})\sqrt{\pi}   \big(\frac{2}{|\m{x}|}\big)^\frac{p-2}{2} I_{\frac{p}{2}}(|\m{x}|)\\
&=\sqrt{\pi}\lambda_{p-1}2^\frac{p-2}{2}\Gamma(\frac{p-1}{2}) \frac{e^{i\langle \m{y},\m{s}\rangle}  }{|\m{x}|^\frac{p}{2} } I_{\frac{p}{2}}(|\m{x}|)\m{x}.
\end{align*}
\begin{flushright}
$\square$
\end{flushright}
\textbf{Conclutions}
In this paper we define a  new subclass for the class of monogenic functions in the biaxial case. This function class is called the class of hypermonogenic solutions. After the definition, we provide methods to find hypermonogenic solutions. This is important to be sure that this class of functions is interesting to study. After that, we deduce the Cauchy formula for the upper hemisphere. One interesting topic in future studies would be to find Cauchy formulas for more general manifolds. In the last part, we present a proper definition for hypermonogenic plane wave solutions and deduce two explicit methods to compute these. \\
\\
\textbf{Acknowledgements}
The second author would like to thank all nice people in Krijgslaan for hospitality and patience.

\vspace{1cm}
\textbf{Al\'{i} Guzm\'{a}n Ad\'{a}n}\\
Clifford Research Group \\
Department of Mathematical Analysis \\
Ghent University\\
Krijgslaan 281\\
9000 Gent, Belgium\\
e-mail: {\tt Ali.GuzmanAdan@UGent.be}\\
\\
\textbf{Heikki Orelma}\\
Laboratory of Civil Engineering\\
Tampere University of Technology\\
Korkeakoulunkatu 10, \\
33720 Tampere, Finland\\
e-mail: {\tt Heikki.Orelma@tut.fi}\\
\\
\textbf{Franciscus Sommen}\\
Clifford Research Group\\
Department of Mathematical Analysis\\
Ghent University\\
Krijgslaan 281\\
9000 Gent, Belgium\\
e-mail: {\tt Franciscus.Sommen@UGent.be}

\end{document}